\documentclass[12pt,onecolumn ]{IEEEtran}

\IEEEoverridecommandlockouts
 \usepackage{geometry}
\usepackage{amsmath, amssymb, amsfonts}
\usepackage{enumerate}
\usepackage{setspace}
\usepackage{graphicx}
\usepackage{color}
\usepackage{verbatim}
\usepackage{url}
\usepackage{algorithm}
\usepackage[noend]{algpseudocode}
\usepackage{float}
\usepackage{amsthm,thmtools}

\makeatletter
\def\BState{\State\hskip-\ALG@thistlm}
\makeatother

\def\urltilde{\kern -.15em\lower .7ex\hbox{\~{}}\kern .04em}
\def\urldot{\kern -.10em.\kern -.10em}
\def\urlhttp{http\kern -.10em\lower -.1ex\hbox{:}\kern -.12em\lower 0ex\hbox{/}\kern -.18em\lower 0ex\hbox{/}}

\declaretheorem[name={Definition}  ] {definition} 
\declaretheorem[name={Theorem}  ] {theorem}

\declaretheorem[name={Problem}  ] {problem} 
\declaretheorem[name={Claim}  ] {claim}

\declaretheorem[name={Corollary}  ] {corollary}

\declaretheorem[name={Remark}  ] {remark}
\declaretheorem[name={Fact}  ] {fact} 

\declaretheorem[name={Proposition}  ] {proposition}

\declaretheorem[name={Example},qed={\lower-0.3ex\hbox{$\square$}} ] {example}

\newcommand{\N}{\mathcal N}

\newcommand{\I}{\mathcal I}

\newcommand{\updt}[1]{{\color{black}#1}}

\begin{document}
\title{A Polynomial-Time Algorithm for Solving the Minimal Observability Problem 
		in Conjunctive Boolean~Networks\thanks{
			This research is supported in part by a research 
		grant from the Israel Science Foundation~(ISF grant 410/15).}} 
		
\author{Eyal Weiss  and Michael Margaliot\thanks{The authors are with the School of Elec.
 Eng., Tel-Aviv University,  Israel 69978. Corresponding Author: michaelm@eng.tau.ac.il}}

\maketitle

 \doublespace

 	 
\begin{abstract}
Many complex systems in biology, physics, and engineering include a large number of state-variables, and measuring the full state of the system is often impossible.  Typically, a set of sensors is used to measure part of the state-variables.
A system   is called observable if these measurements allow to reconstruct the entire state of the system. 
When the system is not observable, 
an important and practical problem is how to add a \emph{minimal} number of sensors so that
the system becomes observable. 
This minimal observability problem is practically useful and theoretically interesting, as it pinpoints
 the most informative nodes in the system.

We consider the minimal observability problem for an important 
 special class of Boolean networks, called  conjunctive Boolean networks~(CBNs).
	Using a graph-theoretic approach,
	we provide a necessary and sufficient condition for observability of a CBN  with~$n$ state-variables, and 
	an efficient~$O(n^2)$-time algorithm for solving the  minimal observability problem. \updt{ We demonstrate the usefulness of these results by studying the properties of a class of random~CBNs.}
\end{abstract}
\section{Introduction}\label{sec:intro}
Real world
 systems   often include    a large number of state-variables~(SVs).  
Measuring all these~SVs   to obtain the  full state of the   system  is often impossible. 
For example, the function of a multipolar neuron 
  may depend upon signals received from thousands of other interconnected neurons
(see, e.g.,~\cite{hawkins2016neurons}).

When the  system has  outputs, that is, functions of~SVs that are directly measured using suitable sensors, 
  an important and  practical question  is whether it is possible 
	to determine the value of \emph{all} the SVs by observing a sequence of the outputs.
	In the language of systems and control theory,   a system that satisfies this property 
	is called \emph{observable}.
	
Establishing observability is the first step in the design of an   \emph{observer}, that is, a 
device that can reconstruct the entire state of the system based on a sequence of the outputs. 
	A typical example is estimating the weather in a certain region based on a small set of
	measurements from local stations in this  region~\cite{farhat2016charney}.
	Observers also  play a crucial  role
in the implementation  of full-state  feedback  controllers  (see, e.g.,~\cite{sontag_book}).

When a  given system is not observable, it is sometimes possible to make it observable by placing
additional sensors that measure
  more (functions of the) SVs. Of course, this may be costly in terms of time and money,
	  so a natural question is:
find the \emph{minimal} number of measurements to add so that the resulting system is observable. 
This \emph{minimal observability problem} is also interesting theoretically, as its
solution means identifying the (functions of) SVs that provide the maximal
 information on the entire state
 of the system~\cite{Liu_observ}.

Minimal observability problems are  recently attracting considerable interest. 
Examples include    monitoring complex services by minimal logging~\cite{biswas2008minimal}, 
the optimal placement of  phasor measurement units   
in power systems (see, e.g.,~\cite{peng2006optimal}), and the minimal sparse 
observability problem addressed in~\cite{sarma2014minimal}. 

Here we solve  the minimal observability problem for an important  special class
of Boolean networks~(BNs). BNs are discrete-time dynamical systems with Boolean SVs. 
BNs have found many applications as models of   dynamical systems. They have been used
to capture the existence and directions of links in 
 complex systems (see, e.g.,~\cite{shahrampour2015topology}),
to  model    social networks (see, e.g.~\cite{soc_nets_bool,soc_nets_plos}), 
 and  the spread of epidemics~\cite{bool_epidemcs}.
In particular, BNs  play  an important role in the modeling of biological processes
and networks, where the feasible set of states is assumed to be finite 
(see, e.g.~\cite{born08,Helikar,faure06}).
A classical example is a gene regulation network,
where each gene may be either ON or OFF (i.e., expressed or not)~\cite{KAUFFMAN1969437}.
In this case, the state of each gene may be   modeled by a 
Boolean~SV, and interactions between the genes 
(through the proteins that they encode) determine
the Boolean update function for each~SV.

Let~$S:=\{0,1\}$. For two integers~$i,j$ let~$[i,j]:=\{i,i+1,\dots,j\}$. 
A BN with $n$ SVs and~$m$ outputs is a discrete-time
dynamical system in the form:
\begin{align}\label{eq:BN}
X_i(k+1)&=f_i(X_1(k),\dots,X_n(k)),  \quad \forall i\in[1,n],\nonumber\\
Y_j(k)&=h_j(X_1(k),\dots,X_n(k)), \quad \forall j\in[1,m]. 
\end{align}
Every~$X_i$ and~$Y_j$ takes values in~$S$,
and $f_i$, $h_j$ are Boolean functions for all~$i,j$,
i.e.,~$f_i, h_j : S^n \to S$.
If there exists an output  
$Y_j(k)=X_i(k)$ then we say that  the SV~$X_i$ 
  is \emph{directly observable}  or  \emph{directly measurable}.

Denote the state of the system at time~$k$
by~$X(k):=\begin{bmatrix} X_1(k)&\dots&X_n(k)\end{bmatrix}'$
 and the output by~$Y(k):=\begin{bmatrix} Y_1(k)&\dots &Y_m(k)\end{bmatrix}'$.
We say that~\eqref{eq:BN} is \emph{observable} on~$[0,N]$,
if any two different initial conditions~$X(0)$ and~$\tilde X(0)$
yield \emph{different} output sequences~$\{ Y(0),\dots,Y(N) \}$ and~$\{\tilde Y(0),\dots, \tilde Y(N) \}$.
This means that 
  given~$\{ Y(0),\dots,Y(N) \}$ it is always possible 
to uniquely determine the initial condition~$X(0)$ of the system.
A BN is called \emph{observable} if it is observable  for some value~$N \geq 0$.

If the output sequence~$\{Y(0),\dots,Y(N)\}$ is identical  for 
two different initial states~$X(0)$ and~$\tilde X(0)$ then it is not possible
to differentiate between them based on the given output sequence.
In this case, we say that these two initial conditions  are \emph{indistinguishable} 
on the time interval~$[0,N]$. 
Clearly, this implies that the BN is not observable on this  time interval.
\updt{
	Boolean control networks (BCNs) are BNs with inputs.
There are several ways to extend the notion of observability to~BCNs
(see, e.g.~\cite{observ_comparison}).
Here, we first consider networks without inputs and then generalize the results 
 to the case with inputs.
}

The observability of   BNs has   been analyzed using 
algebraic and graph-theoretic approaches 
(see, e.g.,~\cite{cheng_book,dima_obs}).
It was proven that testing observability of~BNs is NP-hard 
in the number of~SVs in  the system (see~\cite{dima_obs}). This means that, unless P=NP,
  it is computationally intractable to determine  whether a large~BN is observable.
For a  general survey on the computational complexity of various problems in 
systems and control theory, see~\cite{blondel}.

Of course, the     
  hardness results on determining observability 
	in   general~BNs do  not preclude the possibility that observability
analysis is tractable for 
certain special classes of~BNs.
%
An important    class of~BNs are those with update functions that  include 
nested canalyzing functions~(NCFs) only~\cite{kauffman_93}.
A Boolean function is called canalyzing
if there exists a specific value, called the canalyzing value,
such that an input with this value 
uniquely determines the function's output, regardless of the other variables. 
For example, $0$ is a canalyzing value for the function~AND, 
as~$\text{AND}(0,X_1,\dots,X_k)=0$ for all~$X_i \in \{0,1\}$.
BNs with NCFs are frequently used in
modeling   genetic networks~\cite{harris_2002,kauu_2003,kauu_2004}.  

In this paper, we consider the subclass consisting of NCFs 
which are constructed exclusively by AND operators (i.e. by conjunctive functions).
As models for       gene regulation networks,  conjunctive functions  encode 
\emph{synergistic} regulation of genes by transcription factors~\cite{Jarrah2010}, 
and it seems that this  mechanism indeed exists in certain
regulatory networks~\cite{Merika01051995,Gummow2006,Nguyen2006.0012}.

A BN is called a \emph{conjunctive Boolean network}~(CBN) 
if every update function includes   AND operations only, 
i.e., the dynamics is:
\begin{equation}\label{eq:CBN}
X_i(k+1)=\prod_{j=1}^{n}(X_{j}(k))^{\epsilon_{ji}},\quad \forall i\in [1,n], 
\end{equation}
where~$\epsilon_{ji}\in \{0,1\}$ for all~$i,j$. 
The special (but not very interesting) case~$X_i(k+1)=X_i(k)$ is called 
  a constant updating function.

\begin{remark}\label{rem:orgate}
Note that a BN is called 
a disjunctive Boolean network~(DBN) if every
update function includes only OR operators.
By applying De Morgan Law's, it is possible to reduce~DBNs to~CBNs,
and therefore all the results  in this paper hold for~DBNs as well.
\end{remark}

A useful representation of a~CBN is given by a \emph{dependency graph} (also known as the  wiring diagram).
  This is a directed graph (digraph) 
in which every vertex corresponds to an~SV of the CBN,
and a directed edge~$(i \to j)$ exists if~$X_i(k)$ is one of the arguments 
in  the update function of~$X_j(k+1)$. 
Thus, the dependency graph encodes the variable dependencies in the update functions.

There is a one-to-one correspondence between a~CBN and its dependency graph,
which enables a graph-theoretic analysis of~CBNs. This has been used to
 analyze various 
properties of~CBNs including:
\updt{
	characterization of the periodic orbits~\cite{Jarrah2010,weakly_con_cbn}, 
	robustness of these orbits to single bit perturbations~\cite{basar_CBN},
	and  controllability of~CBNs~\cite{min_control,basar_min_cont}. 
}

However, observability  and, in particular, the minimal observability problem
in  CBNs  has not   been studied before. 
We consider the following problem.

\begin{problem}\label{prob:minobs}
	Given a CBN with~$n$ SVs determine a \emph{minimal} set of indices~$\I \subseteq [1,n]$,
	such that making each~$X_i(k)$, $i \in \I$, 
	directly measurable  yields an observable~CBN.
\end{problem}
Note that an efficient  (i.e.,   polynomial time) solution to Problem~\ref{prob:minobs}
must entail an efficient  algorithm for testing observability of a~CBN.

\begin{example}
	Consider Problem~\ref{prob:minobs} for the~CBN:
	\begin{align*}
	X_1(k+1)&=X_2(k),\\
	X_2(k+1)&=X_1(k)X_2(k). 
	\end{align*}
	Suppose that we make~$X_1(k)$ directly measurable,   that is, add an  output~$Y_1(k)=X_1(k)$.  
	Then the resulting one-output CBN is observable. Indeed, given~$\{Y_1(0), Y_1(1)\}$,
	 the initial condition of the CBN 
	  is~$X_1(0)=Y_1(0)$,~$X_2(0)=Y_1(1)$. Since observability
	requires at least one output, it is clear  that 
	this is a minimal solution to Problem~\ref{prob:minobs}.
\end{example}

\updt{
To make things more concrete consider the following application. 
  A graph describes a  network of interacting  agents 
 with directed edges
 describing  the neighboring relations.
  Every  agent  has two possible opinions on some  matter.
At time~$k$ the opinion of agent~$i$ is described by the state-variable~$x_i(k)\in\{0,1\}$. 
Every agent is ``conservative'' in the sense that it
 tends to hold the opinion zero, unless \emph{all} 
its neighbors hold the opinion one
at time~$k$ and then he updates his opinion to~$x_i(k+1)=1$. 
 Initially, there are no observation nodes, but it is possible to recruit agents
so that they provide reports  on their opinion at any time~$k$. However, 
the recruitment of an  agent is costly in terms  of money, time, etc.
Then a natural  question is: what is the \emph{minimal}
number of agents that must be recruited in order to be able to infer, using a time sequence
of their reports,
 the entire state of the network? This is exactly Problem~\ref{prob:minobs}.

As another application, based on  Remark~\ref{rem:orgate}, consider a model of epidemics
that includes a set of agents that can be either susceptible or infected. 
The directed dependency network  describes contacts between  agents that can lead to infection. The infection is 
so contagious   that a susceptible becomes infected   if even a  single neighbor  is  infected. 
In this context  Problem~\ref{prob:minobs} again has a natural interpretation. 
}

The   contributions of this paper are: 
	\begin{enumerate}
	\item
  a necessary and sufficient condition for the observability of a CBN;
  \item a procedure for designing an observer for an observable CBN;
	and 
	\item
	an~$O(n^2)$-time algorithm for solving Problem~\ref{prob:minobs}.
\end{enumerate}

\updt{
The remainder of this paper   is organized as follows. 
Section~\ref{sec:pre} reviews some standard definitions and notations 
from graph theory that will be used later on. 
Section~\ref{sec:main} describes our main theoretical results.
As already noted by Kauffman~\cite{KAUFFMAN1969437}, there are good reasons to model 
various biological processes using networks of randomly connected binary devices. 
In Section~\ref{sec:random_cbn}, we    use our algorithm to 
 solve Problem~\ref{prob:minobs}
for a class of random~CBNs. These are described by random
dependency graphs with equiprobable edges. 
Surprisingly, perhaps, we show that to make 
 these~CBNs observable, one must observe
 at least~$69\%$ of the nodes.  
Section~\ref{sec:exte} depicts two extensions of our results.
Section~\ref{sec:conclusion} concludes and presents directions for further research.
A detailed description of the main algorithm introduced in the paper 
is given in the Appendix.
}

\section{Preliminaries}\label{sec:pre}
Let $G=(V,E)$ be a digraph, with~$V$ the set of vertices, 
and $E$ the set of directed edges (arcs). Let~$e_{i\to j}$ (or $(v_i\to v_j)$) denote
the arc from~$v_i$ to~$v_j$. When such an arc exists, we say that~$v_i$ is an \emph{in-neighbor} of~$v_j$, and~$v_j$ as an \emph{out-neighbor} of~$v_i$. 
The set of in-neighbors [out-neighbors] of~$v_i$ is denoted by~$\N_{in}(v_i)$ [$\N_{out}(v_i)$]. The \emph{in-degree} [\emph{out-degree}] 
of~$v_i$ is~$|\N_{in}(v_i)|$ [$|\N_{out}(v_i)|$].
A \emph{source} [\emph{sink}] is a node with  in-degree [out-degree] zero.

Let $v_i$ and $v_j$ be two vertices in $V$. A \emph{walk} from~$v_i$ to~$v_j$, denoted~$w_{ij}$, is a sequence: 
$v_{i_0}v_{i_1}\dots v_{i_q}$, with $v_{i_0}=v_i$, $v_{i_q}=v_j$,  
and~$e_{i_k\to i_{k+1}} \in E $ for all $k\in [0,q-1]$. 
A \emph{simple path} is a walk with pairwise distinct vertices. 
We say that $v_i$ is \emph{reachable} from $v_j$ if there exists a simple path from~$v_j$   to~$v_i$. 
A \emph{closed walk} is a walk that starts and terminates  at the same vertex. 
A closed walk is called a \emph{cycle} if all the vertices in the walk are distinct, except for the start-vertex and the end-vertex.

Given a CBN in the form~\eqref{eq:CBN}, the associated \emph{dependency graph} 
is a digraph $G=(V,E)$ with~$n$ vertices (corresponding to the SVs of the system), 
such that~$e_{i\to j} \in E$ if and only if (iff)~$\epsilon_{ij}=1$.
A node   in the dependency graph that
 represents a [non] directly observable SV    
  is called a \emph{[non] directly observable node}.  

\section{Main Results}\label{sec:main}
From hereon, we
  consider CBNs with~$n$ SVs and~$m\geq 0 $ outputs:
\begin{align}\label{eq:CBNee}
X_i(k+1)&=f_i(X_1(k),\dots,X_n(k)),   & \forall i \in [1,n],\nonumber\\
Y_j(k)&=X_j(k),   &\forall j \in [1,m],  
\end{align}
where the~$f_i$s are AND operators,
 and every output~$Y_i$  is the value of an SV (without loss of generality, we   assume that the~$m$ outputs
correspond to the first~$m$ SVs). Thus, nodes~$X_1,\dots,X_m$ [$X_{m+1},\dots,X_n$] in the dependency graph
are   [non] directly observable.

We begin by  deriving two simple   necessary 
conditions for observability of~\eqref{eq:CBNee}.

\begin{definition} 
\updt{
We say that a CBN   has \emph{Property~$O_1$}
if for every non-directly observable node~$X_i$ 
there exists some other node~$X_j$, such that
$\N_{in}(X_j)=\{X_i\}$.
}
\end{definition}

\begin{fact}\label{fact:necessary_1}
	If a CBN is 	  observable then it  has Property~$O_1$.
\end{fact}

{\sl Proof of Fact~\ref{fact:necessary_1}.}
Consider  a CBN that does not satisfy Property~$O_1$. Then it admits 
a non-directly observable node~$X_i$ in its dependency graph,
that  is  \emph{not}    
the only element in the in-neighbors' set of some other   node.
This implies one of the following two cases.

\noindent Case 1: The node~$X_i$ is a sink. Then clearly the~CBN is not observable, 
as there is no way to determine~$X_i(0)$. 

\noindent Case 2: There exists some other node~$X_j$ such that~$N_{in}(X_j)$ 
contains~$X_i$ and at least one  other  node. 
Consider two initial conditions: one with all~SVs  equal to zero, and the second
with all~SVs equal to zero, except for~$X_i(0)$ that is one.
Then for both these conditions
the value of every directly observable node will be zero for all time~$k$, so  these two  
 states are indistinguishable.~\IEEEQED

\begin{figure}[t]
	 	\centering
	
	\begingroup%
	\makeatletter%
	\providecommand\color[2][]{%
		\errmessage{(Inkscape) Color is used for the text in Inkscape, but the package 'color.sty' is not loaded}%
		\renewcommand\color[2][]{}%
	}%
	\providecommand\transparent[1]{%
		\errmessage{(Inkscape) Transparency is used (non-zero) for the text in Inkscape, but the package 'transparent.sty' is not loaded}%
		\renewcommand\transparent[1]{}%
	}%
	\providecommand\rotatebox[2]{#2}%
	\ifx\svgwidth\undefined%
	\setlength{\unitlength}{102.54215871bp}%
	\ifx\svgscale\undefined%
	\relax%
	\else%
	\setlength{\unitlength}{\unitlength * \real{\svgscale}}%
	\fi%
	\else%
	\setlength{\unitlength}{\svgwidth}%
	\fi%
	\global\let\svgwidth\undefined%
	\global\let\svgscale\undefined%
	\makeatother%
	\begin{picture}(1,0.67387395)%
	\put(0,0){\includegraphics[width=\unitlength]{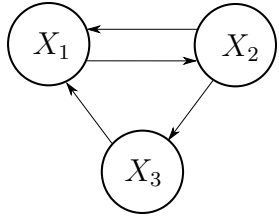}}%
	\put(-0.2712273,0.90792654){\color[rgb]{0,0,0}\makebox(0,0)[lt]{\begin{minipage}{0.47018981\unitlength}\raggedright \end{minipage}}}%
	\put(-1.24295302,2.76778866){\color[rgb]{0,0,0}\makebox(0,0)[lt]{\begin{minipage}{1.25383975\unitlength}\raggedright \end{minipage}}}%
	\put(0.09541399,0.5989266){\color[rgb]{0,0,0}\makebox(0,0)[lb]{\smash{$X_1$}}}%
	\put(-0.02045932,2.4125341){\color[rgb]{0,0,0}\makebox(0,0)[lt]{\begin{minipage}{0.06269197\unitlength}\raggedright \end{minipage}}}%
	\put(0.83633098,-0.46084841){\color[rgb]{0,0,0}\makebox(0,0)[lt]{\begin{minipage}{1.81806762\unitlength}\raggedright \end{minipage}}}%
	\put(1.04530439,-0.42950239){\color[rgb]{0,0,0}\makebox(0,0)[lb]{\smash{}}}%
	\put(0.78846212,0.59249677){\color[rgb]{0,0,0}\makebox(0,0)[lb]{\smash{$X_2$}}}%
	\put(0.43164687,0.1343235){\color[rgb]{0,0,0}\makebox(0,0)[lb]{\smash{$X_3$}}}%
	\end{picture}%
	\endgroup%

	\caption{Dependency graph of the CBN  in Example~\ref{ex:o1}.}
	\label{fig:example1}
\end{figure}

\begin{example}\label{ex:o1}
	Consider the   CBN:
	\begin{align}\label{eq:ftyh}
		X_1(k+1)&=X_2(k)X_3(k),\nonumber \\
		X_2(k+1)&=X_1(k),\\
		X_3(k+1)&=X_2(k),\nonumber\\
		Y_1(k) &=X_1(k)\nonumber.
	\end{align}
	The dependency graph of this CBN  
	does not satisfy Property~$O_1$ (see Fig.~\ref{fig:example1}). 
	Indeed,~$X_3$, which  is a 
	non-directly observable node,   is not 
	the only element in the in-neighbors  set of some other node.
	It is clear that the two initial conditions~$X(0)=\begin{bmatrix} 0& 0& 0 \end{bmatrix}'$
	and~$X(0)=\begin{bmatrix} 0&0&1 \end{bmatrix}'$   are indistinguishable, as for both 
	conditions the output is~$Y_1(k)=0$ for all~$k\geq 0$. 
\end{example}

\begin{definition}
We say that a CBN   has  \emph{Property~$O_2$} if every cycle~$C$ in its dependency graph 
 that  is composed solely   of non-directly observable nodes  satisfies the following property.
$C$ includes  a node~$X_i$ which is the only element in the in-neighbors  set of some other
node~$X_j$, i.e.~$N_{in}(X_j)=\{X_i\}$, and~$X_j$ is not part of the cycle~$C$.
\end{definition}

\begin{fact}\label{fact:necessary_2}
	If a CBN is 	  observable then it  satisfies  Property~$O_2$.
\end{fact}

{\sl Proof of Fact~\ref{fact:necessary_2}.}
\updt{
	Consider a CBN that does not satisfy Property~$O_2$. 
	Then its dependency graph admits 
	a cycle~$C$, composed  solely of non-directly observable nodes,
	and none of these nodes is  
	the only element in the in-neighbors  set of a node that is not part of the cycle~$C$.
	Consider two initial conditions. One with all SVs equal to zero.
	The second with all~SVs equal to zero, 
	except for one~SV that belongs to~$C$, that is equal to one.
	Then these two initial conditions are indistinguishable.~\IEEEQED
}

\begin{figure}[t]
	 	\centering
	  \begingroup%
  \makeatletter%
  \providecommand\color[2][]{%
    \errmessage{(Inkscape) Color is used for the text in Inkscape, but the package 'color.sty' is not loaded}%
    \renewcommand\color[2][]{}%
  }%
  \providecommand\transparent[1]{%
    \errmessage{(Inkscape) Transparency is used (non-zero) for the text in Inkscape, but the package 'transparent.sty' is not loaded}%
    \renewcommand\transparent[1]{}%
  }%
  \providecommand\rotatebox[2]{#2}%
  \ifx\svgwidth\undefined%
    \setlength{\unitlength}{186.86645508bp}%
    \ifx\svgscale\undefined%
      \relax%
    \else%
      \setlength{\unitlength}{\unitlength * \real{\svgscale}}%
    \fi%
  \else%
    \setlength{\unitlength}{\svgwidth}%
  \fi%
  \global\let\svgwidth\undefined%
  \global\let\svgscale\undefined%
  \makeatother%
  \begin{picture}(1,0.35840483)%
    \put(0,0){\includegraphics[width=\unitlength]{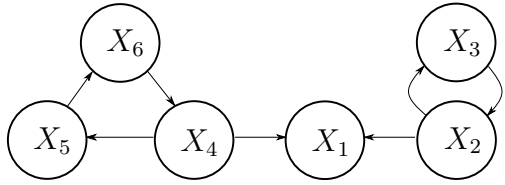}}%
    \put(0.52294758,-0.65171694){\color[rgb]{0,0,0}\makebox(0,0)[lt]{\begin{minipage}{0.25801463\unitlength}\raggedright \end{minipage}}}%
    \put(-0.01028271,0.36887436){\color[rgb]{0,0,0}\makebox(0,0)[lt]{\begin{minipage}{0.68803914\unitlength}\raggedright \end{minipage}}}%
    \put(0.66055542,0.17392995){\color[rgb]{0,0,0}\makebox(0,0)[lt]{\begin{minipage}{0.03440195\unitlength}\raggedright \end{minipage}}}%
    \put(1.13071539,-1.40282626){\color[rgb]{0,0,0}\makebox(0,0)[lt]{\begin{minipage}{0.99765674\unitlength}\raggedright \end{minipage}}}%
    \put(1.24538855,-1.38562525){\color[rgb]{0,0,0}\makebox(0,0)[lb]{\smash{}}}%
    \put(0.05318047,0.06761505){\color[rgb]{0,0,0}\makebox(0,0)[lb]{\smash{$X_5$}}}%
    \put(0.34950954,0.06425774){\color[rgb]{0,0,0}\makebox(0,0)[lb]{\smash{$X_4$}}}%
    \put(0.20043406,0.26493551){\color[rgb]{0,0,0}\makebox(0,0)[lb]{\smash{$X_6$}}}%
    \put(0.61469143,0.06425763){\color[rgb]{0,0,0}\makebox(0,0)[lb]{\smash{$X_1$}}}%
    \put(0.88130634,0.06855773){\color[rgb]{0,0,0}\makebox(0,0)[lb]{\smash{$X_2$}}}%
    \put(0.87987284,0.26493573){\color[rgb]{0,0,0}\makebox(0,0)[lb]{\smash{$X_3$}}}%
  \end{picture}%
\endgroup%
	\caption{Dependency graph of the CBN in Example~\ref{ex:o2}.}
	\label{fig:example2}
\end{figure}

\begin{example}\label{ex:o2}
	Consider the   CBN:
	\begin{align}\label{eq:dyncbn}
		X_1(k+1)&=X_2(k)X_4(k),\nonumber\\
		X_2(k+1)&=X_3(k),\nonumber\\
		X_3(k+1)&=X_2(k),\\
		X_4(k+1)&=X_6(k),\nonumber \\
		X_5(k+1)&=X_4(k),\nonumber\\
		X_6(k+1)&=X_5(k),\nonumber\\
		Y_1(k) &=X_1(k)\nonumber. 
	\end{align}
	This  CBN     has Property~$O_1$, but the cycle formed of 
	$X_4,X_5,X_6$ implies that 
	it does not satisfy  Property~$O_2$ (see Fig.~\ref{fig:example2}). 
	Here the two initial conditions
$
		\begin{bmatrix} 0 &0 &0& 0& 0& 0 \end{bmatrix}'$, and~$\begin{bmatrix} 0 &0 &0& 1& 1&  1\end{bmatrix}'$ 
	yield the same output sequence, namely,~$Y_1(k)=0$ for all~$k\geq 0$, so this CBN is not observable.  
\end{example}

\subsection{Necessary and Sufficient Condition for Observability  }\label{ssec:conditions}
Facts~\ref{fact:necessary_1} and~\ref{fact:necessary_2} provide two
necessary conditions for   observability of a~CBN. 
The next result shows that the combination of these
conditions provides a necessary and sufficient condition for observability.

\begin{theorem}\label{theorem:obs_conditions}
	A CBN is observable iff    it satisfies Properties~$O_1$ and~$O_2$.
\end{theorem}

To prove this, we introduce another definition and several  auxiliary results. 
An \emph{observed  path} in the dependency graph is a non-empty  ordered   set of nodes  
such that: (1)~the last element in the set is a directly observable node; and (2)~if 
the set contains~$p>1$ elements, then for any~$i<p$ 
the~$i$-th element is a non-directly observable node, 
and is the \emph{only} element in the in-neighbors  set of node~$i+1$.
Roughly speaking, an observed path corresponds to a shift register
whose last cell is directly observable. 
Observed paths with non-overlapping nodes are called \emph{disjoint observed paths}.

\begin{proposition}\label{prop:vertex_cover_by_OPs}
	Consider a CBN with a dependency graph~$G$
	that satisfies Properties~$O_1$ and~$O_2$.
	Then~$G$ can be decomposed  
	into disjoint observed paths, such that every vertex in the graph  
	belongs to a single observed path 
	(i.e., the union of the disjoint observed paths is a vertex cover of~$G$).
\end{proposition}

{\sl Proof of Prop.~\ref{prop:vertex_cover_by_OPs}}.
We give a constructive proof.  Algorithm~1   below 
accepts such a graph~$G$ and terminates after each vertex in the graph 
belongs to exactly one observed path.

\begin{algorithm}[H] \label{alg:one} 
\caption{Decompose  the nodes of~$G$ into disjoint observed paths} 
\renewcommand{\algorithmicrequire}{\textbf{Input:}}
\renewcommand{\algorithmicensure}{\textbf{Output:}}
\begin{algorithmic}[1]
\Require Dependency graph~$G$ of a CBN in the form~\eqref{eq:CBNee}
 that satisfies Properties~$O_1$ and~$O_2$.
\Ensure A decomposition of~$G$ into $m$ disjoint observed paths.
\For  {$i=1$ to $m$}
  * every iteration  builds a new path ending with~$X_i$ * 
\State  $\textit{o-node} \gets X_i$ ;
   $\textit{o-path}\gets \{X_i\}$
\If {$|\N_{in}(\textit{o-node})|=1$} \label{if:in_deg_equals_one}
\State Let $v$ be   such that~$\{v\}=\N_{in}(\textit{o-node})$
	\If {$v  $ does not belong to a  previous path  and is not 
		 directly observable}  
		\State   insert~$v$ to $\textit{o-path}$ just before~$\textit{o-node}$
		\State $\textit{o-node} \gets v$;
				  goto~\ref{if:in_deg_equals_one}
	\Else { print} $\textit{o-path}$
	\EndIf 
\EndIf
\EndFor
\State  {\bf end for}
\end{algorithmic}
\end{algorithm}

\updt{We now prove the correctness of this algorithm. 
To simplify the notation, let us say that~$X_p$ \emph{points to}~$X_q$ if~$p\not =q$ and~$N_{in}(X_q)={X_p}$, and
  denote this by~$X_p \mapsto X_q$. The special arrow indictes that the dependency graph includes an edge from~$X_p$ to~$X_q$
	and that there are no other edges pointing to~$X_q$. 

If all the nodes are directly observable (i.e. if~$m=n$) 
  the algorithm will assign every node to a different observed path and this is correct.
Thus, we may assume that~$m<n$. 

Pick   a non directly observable node~$X_j$. Then~$m<j\leq n$. 
Our first goal is to  prove the  following result.
\begin{claim}\label{cla:out}
 The algorithm outputs  an observed path  that contains~$X_j$. 
\end{claim}
By Property~$O_1$, there exists~$k\not = j$ such that  $X_j \mapsto X_k$. 
  We consider two cases. 

\noindent {\sl Case 1.} 
If~$k\leq m$ then~$X_k$ is directly observable and
  the algorithm will add~$X_j$ to an observed path as  
	it ``traces back'' from~$X_k$ unless~$X_j$ has already been included in some other observed path found by the algorithm. 
	Thus, in this case Claim~\ref{cla:out} holds. 

\noindent {\sl Case 2.}  Suppose that~$k>m$, i.e.~$X_k$ is non directly observable.
By Property~$O_1$, there exists~$h\not = k$ such that~$ X_k \mapsto X_h $, so~$X_j \mapsto X_k  \mapsto X_h $.  
If~$h\leq m$ then we conclude as in Case~1 that 
 the algorithm outputs  an observed path  that contains~$X_j$. 
Thus, we only need to consider the case where as we proceed from~$X_j$ using Property~$O_1$ we never ``find'' a
directly observable node. Then 
there exists a set of non directly observable nodes~$X_{k_1},\dots, X_{k_\ell}$, with~$k_1=j$,
 such that 
\begin{align*}
													  X_{k_1} \mapsto X_{k_2} \mapsto \dots \mapsto 
														X_{k_\ell} 	\mapsto	   X_{k_1}  .
\end{align*}
This means that~$X_j$ is   part of a cycle~$C$ of non directly observable nodes. 
 By Property~$O_2$, $C$ includes  a node~$X_{k_{i}}$ such that~$ X_{k_{i}} \mapsto X_{s_1}  $,
  where~$X_{s_1}$ is not part of the cycle~$C$.
If~$X_{s_1}$ is   directly observable then we conclude that the algorithm will output 
an   observed path that includes~$X_j$. 
If~$X_{s_1}$ is not directly observable then by  Property~$O_1$,  there exists~$s_2 \not =s_1 $ such that
$   X_{s_1} \mapsto X_{s_2}$. Furthermore, since every node in~$C$ has in degree one,~$X_{s_2} \not \in C $. Proceeding this way, 
we conclude that there exist~$s_1,\dots , s_p$ such that
\[
					X_{k_{i}} \mapsto X_{s_1}  \mapsto X_{s_2}\mapsto  \dots\mapsto X_{s_p},
\]
with~$X_{s_p}$ a directly observable node. This means that the algorithm will output~$X_j$
in an observed path  as it traces back from~$X_{s_p}$, unless it already included~$X_j$  
 in another observed path. 
This completes the proof of  Claim~\ref{cla:out}.~\IEEEQED

Summarizing, we showed that  \emph{every} non directly observable node~$X_j$ is contained in an observed path produced by the algorithm. 
The fact that every~$X_j$ will be in a single observed path, and that the observed paths will be distinct is clear from the description of the algorithm. 
The algorithm's correctness completes the proof of Prop.~\ref{prop:vertex_cover_by_OPs}.~\IEEEQED  
}

 We can now prove Thm.~\ref{theorem:obs_conditions}. 
 
{\sl Proof of Thm.~\ref{theorem:obs_conditions}.}
\updt{
Consider the following set of statements.
\begin{enumerate}[(a)]
\item The CBN is observable;
\item The dependency graph has Properties~$O_1$ and~$O_2$;
\item There   exists a decomposition 
of the dependency graph  
into a set of~$m\geq 1$ disjoint observed paths~$O^1,\dots,O^m$, 
such that every vertex in the graph belongs to a single observed path.
\end{enumerate}

We already know that~(a) $\to$ (b). The   correctness of Algorithm~1 
implies that  (b) $\to$ (c). If~(c) holds
then  the values of the output of~$O^i$
at times~$0, \dots, N_i-1$ are the initial values of the SVs in~$O^i$,
organized in reverse  order. 
Therefore  it is possible to determine the initial condition of every~SV
in  the~CBN using the output sequence on~$[0,\max_{i=1,\dots,m} \{N_i \} -1]$.
Thus, the~CBN is observable, so (c) $\to$ (a). We conclude that statements (a), (b), and~(c) are all equivalent and this proves
Thm.~\ref{theorem:obs_conditions}.~\IEEEQED

} 

The proof of Thm.~\ref{theorem:obs_conditions}    implies the following. 
 \begin{corollary}
	A CBN is observable iff  its dependency graph can be decomposed  
	into a set of disjoint observed paths.
\end{corollary}

The proof of Thm.~\ref{theorem:obs_conditions}     also provides a way to  design 
an observer 
for  an observable CBN. The procedure is as follows:
\begin{enumerate} [(a)]
	\item Construct the dependency graph~$G$;
	\item Apply Algorithm~$1$ to decompose the nodes of~$G$  into a set of disjoint observed paths;
	\item Observe an output sequence of length equal to the   longest observed path;
	\item Map   the values observed at each output to the values of the SVs composing the observed paths, 
	in reverse order, to obtain the initial state of the entire CBN.
\end{enumerate}

\noindent Of course, once the initial condition is recovered,
the known dynamics of the~CBN allows to determine the state of the~CBN at any time step. 

\updt{

\begin{example}\label{ex:alg_one}
	Consider the CBN:
	\begin{align}
	X_1(k+1)&=X_3(k),\nonumber \\
	X_2(k+1)&=X_5(k),\nonumber\\
	X_3(k+1)&=X_4(k),\nonumber\\
	X_4(k+1)&=X_2(k)X_3(k),\\
	X_5(k+1)&=X_1(k)X_5(k),\nonumber\\
	Y_1(k) &=X_1(k),\nonumber\\
	Y_2(k) &=X_2(k)\nonumber.
	\end{align}
	The dependency graph of this CBN  
	satisfies Properties~$O_1,O_2$ (see Fig.~\ref{fig:alg_one}),
	so Thm.~\ref{theorem:obs_conditions} implies that it is observable,
	and decomposable to a set of disjoint observed paths.
	Applying Algorithm~$1$ to this CBN yields~$O^1=(X_4,X_3,X_1)$,
	$O^2=(X_5,X_2)$, where~$X_4 \mapsto X_3 \mapsto X_1$,
	$X_5 \mapsto X_2$. 
\end{example}

\begin{figure}[t]
	\centering
	
	\begingroup%
	\makeatletter%
	\providecommand\color[2][]{%
		\errmessage{(Inkscape) Color is used for the text in Inkscape, but the package 'color.sty' is not loaded}%
		\renewcommand\color[2][]{}%
	}%
	\providecommand\transparent[1]{%
		\errmessage{(Inkscape) Transparency is used (non-zero) for the text in Inkscape, but the package 'transparent.sty' is not loaded}%
		\renewcommand\transparent[1]{}%
	}%
	\providecommand\rotatebox[2]{#2}%
	\ifx\svgwidth\undefined%
	\setlength{\unitlength}{115.56282824bp}%
	\ifx\svgscale\undefined%
	\relax%
	\else%
	\setlength{\unitlength}{\unitlength * \real{\svgscale}}%
	\fi%
	\else%
	\setlength{\unitlength}{\svgwidth}%
	\fi%
	\global\let\svgwidth\undefined%
	\global\let\svgscale\undefined%
	\makeatother%
	\begin{picture}(1,1.30818466)%
	\put(0,0){\includegraphics[width=\unitlength]{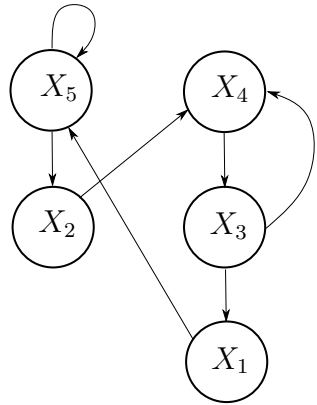}}%
	\put(-0.30087243,1.38377171){\color[rgb]{0,0,0}\makebox(0,0)[lt]{\begin{minipage}{0.41721269\unitlength}\raggedright \end{minipage}}}%
	\put(-1.16311209,3.03407986){\color[rgb]{0,0,0}\makebox(0,0)[lt]{\begin{minipage}{1.11256739\unitlength}\raggedright \end{minipage}}}%
	\put(-0.07835893,2.71885246){\color[rgb]{0,0,0}\makebox(0,0)[lt]{\begin{minipage}{0.05562835\unitlength}\raggedright \end{minipage}}}%
	\put(0.68189529,0.16921907){\color[rgb]{0,0,0}\makebox(0,0)[lt]{\begin{minipage}{1.61322271\unitlength}\raggedright \end{minipage}}}%
	\put(0.86732329,0.19703328){\color[rgb]{0,0,0}\makebox(0,0)[lb]{\smash{}}}%
	\put(0.09407721,0.56611273){\color[rgb]{0,0,0}\makebox(0,0)[lb]{\smash{$X_2$}}}%
	\put(0.65643099,1.00631283){\color[rgb]{0,0,0}\makebox(0,0)[lb]{\smash{$X_4$}}}%
	\put(0.09950653,1.0098657){\color[rgb]{0,0,0}\makebox(0,0)[lb]{\smash{$X_5$}}}%
	\put(0.65810603,0.12237715){\color[rgb]{0,0,0}\makebox(0,0)[lb]{\smash{$X_1$}}}%
	\put(0.65267683,0.56117827){\color[rgb]{0,0,0}\makebox(0,0)[lb]{\smash{$X_3$}}}%
	\end{picture}%
	\endgroup%
	
	\caption{Dependency graph of the CBN in Example~\ref{ex:alg_one}.}
	\label{fig:alg_one}
\end{figure}

}

\subsection{Minimal Observability  Problem}\label{ssec:min_obs}

We now use the conditions   in Thm.~\ref{theorem:obs_conditions}  
to efficiently solve   Problem~\ref{prob:minobs}. 
We consider a CBN in the form~\eqref{eq:CBNee},
and the problem is to find a minimal number of additional SVs to measure 
so that the~CBN becomes observable. Of course, if~\eqref{eq:CBNee} is already observable then the solution
to this problem is zero. 

Algorithm~$2$ below solves this   problem. 
For the sake of clarity, we provide here a high-level
  description of the algorithm.
A more detailed description of the algorithm is given in   the Appendix. 

\begin{algorithm}[H] \label{alg:two} 
\caption{Solving the minimal observability problem: a high-level description}  
\renewcommand{\algorithmicrequire}{\textbf{Input:}}
\renewcommand{\algorithmicensure}{\textbf{Output:}}
\begin{algorithmic}[1]
\Require A CBN~\eqref{eq:CBNee} with~$n$ SVs and~$m\geq 0$ outputs.
\Ensure A minimal set of SVs so that
		making these SVs directly observable yields an observable~CBN. 
\State generate the dependency graph~$G=(V,E)$ 
\State create a list~$L_1$   of all~SVs that are 
not directly observable and are not  the only element in the in-neighbors' set of another   node \label{step:genl1}
\State create a list~$L_2$ of all   SVs that are 
not directly observable and   are the  only element in the in-neighbors  set of another   node
\label{step:genl2}
\State create  a list~$L_C$ of   cycles composed solely out of nodes in~$L_2$
\label{step:genlC}
\State  for each cycle~$C\in L_C$, check if one of its elements appears as the
 only element in the in-neighbors  set of another   node that is not part of~$C$.
 If so, remove~$C$ from~$L_C$
 \label{step:reducelC}
\State  copy~$L_1$ into a list~$\I $; pick one  element from each cycle~$C\in L_C$, 
and add these elements to~$\I$   \label{step:pick}
\State return the list~$\I$
\end{algorithmic}
\end{algorithm}

\begin{example}
	Consider the CBN in Example~\ref{ex:o1}. 
	Applying Algorithm~2 to this CBN yields $L_1=\{X_3\}$, $L_2=\{X_2\}$, and~$L_C=\emptyset$.
	The algorithm thus returns~$L_1=\{X_3\}$. Making this a directly observable node yields the CBN
	with dynamics~\eqref{eq:ftyh}
and outputs~$		Y_i(k)=X_i(k)$, $i=1,3$.
	This CBN is indeed observable, and since the algorithm added a single output is is clear that
	this is a minimal solution.
	
	Now consider the CBN in Example~\ref{ex:o2}.
	Applying Algorithm~2 to this CBN yields $L_1=\emptyset$, $L_2=\{X_2,X_3,X_4,X_5,X_6\}$, and~$L_C=\{  \{ X_2,X_3\} , \{ X_4,X_5,X_6\} \}$.
	Thus, Step~\ref{step:pick} in Algorithm~2 yields, say, the output~$ \{X_2,X_4\} $. 
	Making these two nodes directly observable yields the CBN with the dynamics in~\eqref{eq:dyncbn}
	and outputs
	$		Y_i(k) =X_i(k)$, $ i=1,2,4$.
	It is straightforward to verify that this CBN is indeed observable, and also that this addition 
	of two outputs is a solution of the minimal observability problem.  
	
\end{example}

\begin{theorem}\label{theorem:alg_validity}
  Algorithm 2 provides  a   solution to Problem~\ref{prob:minobs}.
\end{theorem}
{\sl Proof of Thm.~\ref{theorem:alg_validity}.}
It is clear that the algorithm always terminates. 
Note that in step~\ref{step:genl1} of the algorithm all the~SVs that do not satisfy
Property~$O_1$ are placed  in the list~$L_1$, and 
in steps~\ref{step:genlC} and~\ref{step:reducelC},
all the cycles that do not satisfy Property~$O_2$ are placed  in~$L_C$, and only those cycles.
Step~\ref{step:pick} initializes~$\I$ as~$L_1$ and then
picks a representative of each cycle in~$L_C$
and then adds it to~$\I$.
Therefore, after making each of the SVs~$X_i(k)$, $i\in \I$,   
directly observable  the modified~CBN
 satisfies the conditions 
in  Thm.~\ref{theorem:obs_conditions},
and hence is observable for  some~$N \geq 0$.

To prove that~$\I$ is   minimal, note 
 that  since~$L_C$ includes only nodes from~$L_2$, it is clear that every cycle in~$L_C$ 
does not include nodes in~$L_1$. Making a node from~$L_1$ directly observable 
 does not change  the fact that
every node in~$L_C$ does not satisfy Property~$O_2$. Therefore, a minimal solution must be as composed by  the algorithm.~\IEEEQED

\subsubsection*{Complexity Analysis of Algorithm~$2$}
Generating the dependency graph~$G$
requires going through~$n$ updating functions, 
and each function has at most~$n$ arguments, so the  complexity of this step is~$O(n^2)$.
The resulting graph satisfies~$|V|=n$, and~$|E| \leq n^2$.
The complexity of each of the other steps in the algorithm is  
  at most linear in~$|V|, |E|$, i.e., it is~$O(n^2)$.
  Summarizing, the complexity of the algorithm is linear in the
length of the description of the~CBN, and the latter is~$O(n^2)$.

Since the algorithm arbitrarily selects  one element from each cycle in~$L_C$, it provides a specific solution to the minimal observability problem.
It is straightforward to modify this so that the algorithm will return   the information needed to build \emph{all} possible solutions. 
Note   that if the algorithm  returns an output list that is empty then the~CBN is observable, so it can also be used to determine if a given~CBN is observable or not. 

\updt{
 \section{Minimal observability in random~CBNs}\label{sec:random_cbn} 

Recall that we can represent a CBN via its dependency graph.
In this section, we consider the case where the  
  dependency graph is  generated 
as a  \emph{directed} Erd\H{o}s-R\'{e}nyi graph~\cite{vanderHofstad:2016:RGC:3086985}, i.e., we fix the number of vertices~$n$  and a probability~$p\in[0,1]$,
and each possible directed edge in the graph is included
with probability~$p$, independently of any other edge.
We then study  the minimal observability problem 
for such random~CBNs via both simulations and analysis.
\subsubsection*{Simulations}
We generated random dependency graphs with~$n=1000$ vertices  
for a set of~$p$ values. 
For every graph we ran the algorithm described here to obtain  the solution~$k$
to the minimal
observability problem, and calculated~$ 100k/n$, i.e. the percentage  of nodes that 
must be added as observed nodes in order to make the~CBN observable. 
For each value of~$p$ we averaged the minimal number of outputs required
over~$100$ independent trials to obtain the average value~$s:=<100k/n>$.
The middle curve in Fig.~\ref{fig:simulation} depicts~$  s $
as a function of~$p$. 

\begin{figure*}[t]
	\centering
	\includegraphics[width=0.75\textwidth]{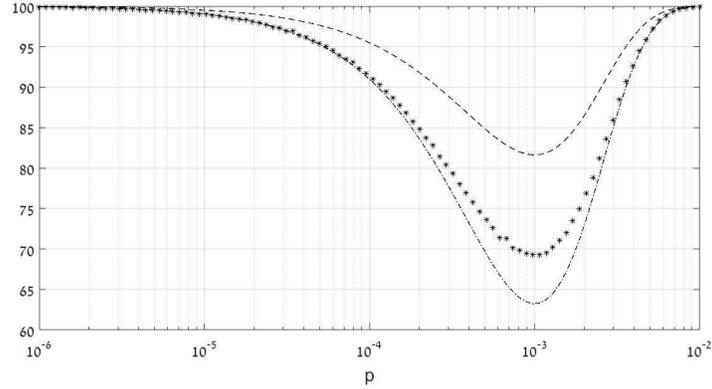}
	\caption{Three functions describing values [in \%] as a function of~$p$.
		Upper curve:~$\bar s (p)$. 
		Middle curve:~$s $ as obtained in the simulations.
		The results for every value~$p$ are based on averaging~$100$ independent random graphs, each with~$1000$ vertices. Lower curve:~$ \underline s(p)$.}
	\label{fig:simulation}
\end{figure*}

It may be seen that~$ s$ decreases sharply around~$p=1/1000 (=1/n)$, 
and achieves a minimum value~$  s^* =69.3\%$, suggesting that the
optimal value for the probability is~$p^*=1/n$.
Similar results were found     when simulating 
for other values of~$n$ in the range~$[10^2,10^4]$. In other words, even in the best
possible random~CBN, on average  about~$70\%$ of the nodes must be added as outputs in order to obtain observability.

The polynomial complexity
of Algorithm~$2$   makes it possible 
to solve  the minimal observability problem even for large values of~$n$.
For example, for a graph with~$n=1000$ the typical running time of the algorithm, 
implemented in MATLAB 
using a standard PC (Intel core i$5$ processor,~$4$GB RAM memory)
is about~$0.03$ seconds. For~$n=10^4$ 
the running time is about~$2.8$ seconds.

We now show how the analytical results described
 in Section~\ref{sec:main} and, in particular, the notion of observed path
allow  to analyze the random CBNs simulated here. 

\subsubsection*{Analysis}
For a dependency graph~$G=(V,E)$, with~$|V|=n$,
let~$k$ denote the minimal number of nodes that must be made directly measurable 
 in order to make the~CBN observable, i.e. the size of the solution to Problem~\ref{prob:minobs}.
We begin with deriving a lower and upper bound on~$k$. 

The key point in the analysis is   the set  of vertices that have in-degree one and no self-loops. Denote this set by~$W$. 
For the random graph described above
 the probability of a node to belong to~$W$ is:
\[
  q(p):=(n-1)p(1-p)^{n-1}, 
\]
so~$E(|W|)=q(p) n $.

The results in Section~\ref{sec:main} 
imply  that every node in~$V\setminus  W$  can   only  be the first node in an  
  observed path.
Since the number of observed paths 
is equal to the number of outputs,~$k\geq n - |W|$.
The exact number of   outputs  needed to achieve observability 
 is determined by the   topology of the graph.
One  optimal case 
is  when  all the vertices of~$W$ form one observed path 
starting with one of the nodes in~$V \setminus W$
yields  the   lower bound:
$k \geq  n - |W|$.\footnote{
One case where this bound is (almost) tight is when~$V = W$, with all the nodes forming a 
 cycle, as in  this case,~$k=1$.}
We conclude that
\begin{equation} \label{eq:lbni}
(100/n) E(k )\geq 100( n - E(|W| ))/n= \underline s,
\end{equation}
with~$\underline s(p):=100(1-q(p))$. The multiplication by~$100/n$ is used to obtain the results 
in terms of percent. 

To derive an upper bound on~$k$, note that every node~$v \in W$ which is located at the 
\emph{beginning} of an observed path  implies a needed  output (in addition to
those necessary for the nodes of~$V \setminus W$).
By the definition of~$W$, its elements have in-degree one.
Therefore, a node~$v \in W$ might be located at the 
beginning of an observed path  only when there is 
a cycle   formed exclusively   of nodes in~$W$.
Let~$C(W)$ denote  the number of cycles composed solely
out of vertices belonging to~$W$. Then we conclude that~$k \leq n - |W| + C(W)$.
Since the nodes of~$W$ do not have self-loops, 
the smallest possible cycle includes two vertices.
Hence,~$C(W) \leq |W|/2$ and this yields~$k \leq n - |W|/2$.
Thus,
\begin{equation}\label{eq:upeer}
(100/n)E(k) \leq 100(n-\frac{E(|W|)}{2})/n=\bar s ,
\end{equation}
where~$\bar s(p):= 100(1-\frac{q(p)}{2}) $. 

Fig.~\ref{fig:simulation} depicts~$\underline s (p) $ and~$\bar s(p)$ as a function of~$p$.
It may be seen that these functions indeed provide a lower and upper bound for
the value~$s$ obtained in the simulations. 

We now turn to determine the optimal value~$p^*$, i.e. the 
edge probability that yields, on average, the smallest solution to the minimal observability problem. 
Clearly,  the probability that maximizes~$q(p)$
minimizes the bounds for~$  s$.
It is straightforward to verify that~$q(p)$ admits a unique maximum at
\[
				p^*:=n^{-1},
\]
and this agrees well with the simulation results.

The simulation results show that for the optimal value~$p=p^*$ 
 the corresponding minimal value is~$  s^*\approx 69.3\%$, 
that is,  even in the optimal case 
about~$70\%$ of the nodes must be observed in order to make 
the~CBN observable. 
To explain  this value, note that
\begin{align*}
				q(p^*) = \left (1- \frac{ 1}{n} \right)^n 
				       \approx e^{-1}.
\end{align*}
Thus, for the  optimal topology  the percentage of outputs needed   on average is:
\begin{align}\label{eq:LB}
	\underline s (p) =100( 1-e^{-1} ) \approx  63.2\%,
\end{align}
 Of course, this is   a lower bound on the number of needed observation nodes, 
as  there is no reason for the optimal   topology to appear frequently in the random simulations.
This analysis agrees well with the simulation results.

The graph of~$s(p)$ in Fig.~\ref{fig:simulation} shows a sharp rise near~$p^*=1/n$. 
To explain this, 
 we compute   the second derivative of~$q(p)$ at~$p=p^*$, namely, 
\begin{align*}
\frac{d^2 q(p) }{dp^2}|_{p=p^*}
&=\left . -(n-1)^2 (1-p)^{n-3} (2-n p) \right |_{p=p^*}\\
&=-(n-1)^2 (1-n^{-1})^{n-3}.
\end{align*}
For~$n \gg 1$ this yields 
\begin{align*}
|\frac{d^2 q(p) }{dp^2}|_{p=p^*} 
&=(n-1)^2 (1-n^{-1})^{n-3} \\ 
&=(n-1)^2 (1-n^{-1})^{ -3} (1-n^{-1})^{n }\\
& \approx  \frac{n^3}{n-1}  e^{-1} .
\end{align*}
This large value of the second derivative implies  a rapid change in the slope of the curve near~$p=p^*$.
Roughly speaking, this means that outside a small 
interval  of probability values around~$1/n$
  a very large number of   outputs is needed to achieve  observability. Again, this agrees well with the simulation results. 
}

\section{Extensions}\label{sec:exte}

We describe two simple extensions of the results above.
\subsection{Observability of CBNs With Inputs}
Consider CBNs with inputs, that is,  
conjunctive Boolean control networks~(CBCNs). 
\updt{
As noted in the introduction, 
there are several possible definitions for observability of BCNs 
with outputs (see, e.g.,~\cite{observ_comparison,cheng_book}).
}
For example, one definition of observability requires that for any two different initial conditions~$a$ and~$b$
 there exists a control sequence
(that may depend on~$a$ and~$b$) guaranteeing that  the output sequences will be different. 
A different possibility is to require that there exists a specific control sequence yielding an output sequence
that distinguishes between any two different initial conditions.

For the case of~CBCNs, it is clear that the ``most informative'' control sequence is when all inputs are one for all time. 
Indeed, any zero input may only obscure the value of an~SV. 
Thus, we use the following definition.
\begin{definition}
A CBCN is said to be observable on~$[0,N]$
if the CBN obtained by setting all inputs to one for all time~$ k\in[0,N-1]$ 
is observable on~$[0,N]$.
%
\end{definition}

This   means that we can first set all inputs~$U_i(k)$ to one, simplify the network by using the fact 
that~$\text{AND}(1,X_1,\dots,X_k)=\text{AND}(X_1,\dots,X_k)$, and then analyze observability and solve the minimal observability 
problem for the resulting~CBN using the approach described in the previous section.

\subsection{Observability of CBNs With More General Output Functions}
Consider a CBN with outputs that are more general than in~\eqref{eq:CBNee}, namely, 
\begin{align}\label{eq:gg}
Y_i(k)&=g_i(X_1(k),\dots,X_n(k)), \quad \forall i\in [1,m],  
\end{align}
with every~$g_i$   an AND operator. 

Consider an augmented BN with~$n+m$ SVs and~$m$ outputs:
\begin{align}\label{eq:cbnaugm}
\bar X_i(k+1)&=f_i(\bar X_1(k),\dots,\bar X_{n}(k)),\quad \forall i \in [1,n],   \nonumber\\
\bar X_{n+j}(k+1)&=  g_j(\bar X_1(k),\dots,\bar X_{n}(k)),  \quad \forall j \in [1,m], \nonumber \\
\bar Y_p(k)&=\bar X_{n+p}(k),\quad \forall p \in [1,m].     
\end{align}
This is a CBN in the form~\eqref{eq:CBNee}, 
 and~$\bar Y_j(k+1)=\bar X_{n+j}(k+1)=g_j(\bar X_1(k),\dots,\bar X_{n}(k))=Y_j(k)$ for all~$j$ and~$k$. 
Thus,~\eqref{eq:cbnaugm}  is observable iff   the CBN with outputs~\eqref{eq:gg} is observable.
In other words, any CBN with outputs in the form~\eqref{eq:gg}, where the~$g_i$s are AND operators, 
can be reduced to the form~\eqref{eq:CBNee}, and then all of the above results 
on observability analysis and minimal observability can be applied.

\section{Discussion}\label{sec:conclusion}
Observability is a fundamental property of dynamical systems, and it  plays a crucial role in the design of
observers, and   full-state feedback controllers. When a system is not observable an important question is 
  to determine a \emph{minimal}  set of measurements so that the system becomes observable. 
	In the context of biological systems, this amounts 
	to determining the minimal number of
	sensors to add so that the measurements will allow to determine the initial state of the biological system. 
	This is important when  the system includes  a large number of SVs and adding sensors is costly in terms 
	of time, money, etc.
 	
	We considered the minimal  observability problem for~CBNs. 
	Using the dependency graph, 
we derived a necessary and sufficient condition for observability of CBNs, and an~$O(n^2)$-time
 algorithm for solving the minimal observability problem for a CBN with~$n$ SVs.
This also includes an explicit procedure which describes the construction 
of an observer for observable CBNs. 

For LTI systems, it is well-known that   controllability analysis and observability analysis are dual problems. 
For nonlinear dynamical
systems, such as BNs, this is not true anymore. Indeed,
it was recently shown that for CBNs the minimal controllability problem   is 
 NP-hard~\cite{min_control} (see also~\cite{basar_min_cont} for some related considerations),
implying that there does not exist an algorithm solving it in polynomial time, 
 unless~P=NP. 

Although the necessary and sufficient conditions for controllability and observability 
of CBNs are quite  analogous (see the definition of a controlled path in~\cite{min_control}), a key difference is that    adding 
a control input to a CBN in order to ``improve'' the controllability 
  changes the dynamics of the CBN 
(and so changes its dependency graph). On the other-hand, adding an output in order to ``improve'' the observability 
does not change the  dynamics.
This is why  the minimal observability problem is computationally  more feasible  than  
   the minimal controllability problem.  

\updt{
The results here suggest several directions for further research. 
Recall that in   \emph{undirected} Erd\H{o}s-R\'{e}nyi graphs
the size of the largest connected component undergoes a phase transition when the edge probability~$p$ crosses the value~$1/n$~\cite[Ch.~4]{vanderHofstad:2016:RGC:3086985}.
Our results show that~$p^*=1/n$ is the ``best'' value  when considering the  minimal observability problem
for CBNs described by directed  Erd\H{o}s-R\'{e}nyi graphs. It may be of interest 
to investigate if~$p^*$ is also the ``best'' value when considering the minimal \emph{controllability} problem
for CBNs, and if other, more general,  BNs demonstrate   some special properties for   this  value of edge probability. 
Another natural direction for future research is the extension of the theoretical results described here 
to more general classes of~BNs.  
  }
	
\section*{Acknowledgment} 
 We are grateful to the anonymous reviewers for many valuable comments.

\section*{Appendix: Detailed Description of Algorithm~2.}

\setcounter{algorithm}{1} 

\begin{algorithm}[H] \label{alg:three} 
\caption{Solving the minimal observability problem }  
\renewcommand{\algorithmicrequire}{\textbf{Input:}}
\renewcommand{\algorithmicensure}{\textbf{Output:}}
\begin{algorithmic}[1]
\Require A CBN~\eqref{eq:CBNee} with~$n$ SVs and~$m\geq 0$ outputs.
\Ensure A minimal set of SVs so that
		making these SVs directly observable yields an observable CBN. 
\State generate the dependency graph~$G=(V,E)$ 
\State initialize   lists~$L_1$ and $L_2$ - each  with~$n$ bits set to zero  ; 
 initialize a matrix~$L_{\textit{pairs}}$ of~$n \times 2$ entries all set to zero
\For {$i=1$ to $n$}  
	  * build $L_2$ and~$L_{\textit{pairs}}$  * \label{step:start_build_l2}
	\If {$|\N_{in}(X_i)|=1$ and $X_j \in \N_{in}(X_i)$ is not directly observable   
		and $i \neq j$} 
		\State $L_2(j) \gets 1$;
	  $L_{\textit{pairs}}(i,1) \gets 1$ ;
		 $L_{\textit{pairs}}(i,2) \gets j$ \label{step:end_build_l2}
	\EndIf	
\EndFor
\For {$i=m+1$ to $n$} \label{step:start_build_l1}
	\Statex * scan over   non-directly observable nodes to build $L_1$  *
	\If {$L_2(i)=0$}
		  $L_1(i) \gets 1$ \label{step:end_build_l1}
	\EndIf
\EndFor
\State copy the list $L_2$ into a list~$L_3$
\For {$i=1$ to $n$}  
	  * build $L_3$  *\label{step:start_build_l3}
	\If {$L_3(i)=0$}
		  $k \gets i$
		\Statex \hspace*{1.2cm}  * $X_k$ is   directly observable   or since it is in~$L_1$ will become     directly observable      * 
		\If {$L_{\textit{pairs}}(k,1)=1$} \label{if:check_sender}
			\State $p \gets L_{\textit{pairs}}(k,2)$
			\State $L_{\textit{pairs}}(k,1) \gets 0$;
			  $L_{\textit{pairs}}(k,2) \gets 0$
			\State $k \gets p$;
			  $L_3(k) \gets 0$  
			\Statex \hspace*{1.55cm} * trace back to    $X_p$ and remove it from $L_3$  *
			\State goto~\ref{if:check_sender} \label{step:end_build_l3}
		\EndIf
	\EndIf
\EndFor
\State generate a digraph~$\tilde{G}$ by removing from~$G$ all the
vertices   that are not in~$L_3$ and all the incident edges  
\State generate a list~$L_C$ of the cycles of~$\tilde{G}$ \label{step:start_division}
\State copy~$L_1$ into a list~$\I$;   pick
 one  element from each cycle in~$L_C$, 
and add       to~$\I$ \label{step:end_division}
\State return the list~$\I$
\end{algorithmic}
\end{algorithm}

The algorithm uses  several  data structures. $L_1$ is a list of~$n$ bits.
For any node~$X_i$ in the dependency graph~$L_1(i)=1$ if~$X_i$ is 
not directly observable and is not the only element in the in-neighbors' set of some other
   node. Otherwise,~$L_1(i)=0$. 
$L_2$ is a list of~$n$ bits with~$L_2(i)=1$ if~$X_i$ is
not directly observable and is the  only element in the in-neighbors' set of some other node. Otherwise,~$L_2(i)=0$. 
$L_{\textit{pairs}}$ is a matrix of dimension~$n \times 2$  such that~$L_{\textit{pairs}}(i,1)=1$ if node~$X_i$ has in-degree one,
i.e.~$N_{in}(X_i)=\{X_j\} $ for some~$j$, 
and in this case~$L_{\textit{pairs}}(i,2)=j$.
The list~$L_3$ includes  vertices that  are part of a cycle,
such that non of the elements composing the cycle appears as the only element in the in-neighbors' set of another node that is not part of the cycle. 
This is created from~$L_2$  using the auxiliary list~$L_{\textit{pairs}}$.


Steps~\ref{step:start_build_l2}-\ref{step:end_build_l2} initialize~$L_2$ 
and~$L_{\textit{pairs}}$.
For convenience, denote the list of directly observable nodes 
by~$L_{\textit{DON}}$ (from the definition of the CBN~\eqref{eq:CBNee}
this  is simply the first~$m$ nodes).
Using~$L_2$ and $L_{\textit{DON}}$, 
  steps~\ref{step:start_build_l1}-\ref{step:end_build_l1} 
form~$L_1$   as
$
L_1 \gets V \setminus \{ L_{\textit{DON}} \cup L_2 \}.
$
Note that the sets~$L_1$, $L_2$, and $L_{\textit{DON}}$ are disjoint, 
with their union equal to the set of vertices of~$G$.
Steps~\ref{step:start_build_l3}-\ref{step:end_build_l3} generate~$L_3$
 by a gradual reduction of~$L_2$
 using the matrix~$L_{\textit{pairs}}$,  
and a depth-first-search-like mechanism.
From $L_3$ it is immediate to obtain the list of cycles~$L_C$. 
Indeed, if~$L_3$ is not empty at the end of
the reduction process, then it is easy to verify that it 
contains exactly the list of vertices composing the desired cycles, 
but not yet divided to sets according to the different cycles.
Steps~\ref{step:start_division}-\ref{step:end_division} 
perform the division to the different sets. This 
  can be implemented   using a strongly connected components algorithm (which is linear in~$|V|,|E|$, see, e.g.,~\cite{Tarjan72depthfirst}), since every 
connected component in the digraph~$\tilde{G}$ is   a cycle.


\section*{Biographies}

{\bf Eyal Weiss}  was born in Israel in 1988. He received the B.Sc. degree (summa cum laude) in
Electrical and Electronic Engineering from Tel Aviv University, Israel, in 2016. He is currently
pursuing his M.Sc. degree at the School of Electrical Engineering-Systems, Tel Aviv University,
Israel.

{\bf Michael Margaliot}  received the BSc (cum laude) and MSc degrees in Elec. Eng. from the
Technion-Israel Institute of Technology-in 1992 and 1995, respectively, and the PhD degree
(summa cum laude) from Tel Aviv University in 1999. He was a post-doctoral fellow in the Dept.
of Theoretical Math. at the Weizmann Institute of Science. In 2000, he joined the Department
of Electrical Engineering-Systems, Tel Aviv University, where he is currently a Professor and
Chair. His research interests include the stability analysis of differential inclusions and switched
systems, optimal control theory, fuzzy control, computation with words, Boolean control networks, contraction theory,
and systems biology. He is co-author of \emph{New Approaches to Fuzzy Modeling and Control: Design and Analysis}, World
Scientific, 2000 and of \emph{Knowledge-Based Neurocomputing}, Springer, 2009. 
He currently serves as an Associate Editor
for the journal \emph{IEEE Transactions on Automatic Control}.

\bibliographystyle{IEEEtran}
\bibliography{eyal_bib}

 \end{document}